\documentclass[12pt, a4paper]{article}
\usepackage[top=2.54cm, bottom= 2.54cm, left=2.8cm, right=2.8cm]{geometry}
\def\Dj{\hbox{D\kern-.73em\raise.30ex\hbox{-}
\raise-.30ex\hbox{}}}
\def\dj{\hbox{d\kern-.33em\raise.80ex\hbox{-}
\raise-.80ex\hbox{\kern-.40em}}}
\usepackage{epsfig}
\usepackage{amsmath,amsthm,amsfonts,amssymb,amscd,cite}

\usepackage{amsfonts,amssymb,color}
\usepackage{mathrsfs}
\usepackage{latexsym}
\usepackage{graphicx}

\allowdisplaybreaks

\newtheorem{thm}{Theorem}
\newtheorem{definition}[thm]{Definition}

\newtheorem{lemma}[thm]{Lemma}

\newtheorem{proposition}[thm]{Proposition}

\begin{document}

\vspace*{20mm}

\noindent {\Large \bf Note on Sombor index of connected graphs \\[2mm]  with given degree sequence}\footnote{This
work is partially supported by the Natural Science Foundation of Guangdong Province (No.~2022A1515011786)
and the National Natural Science Foundation of China (No.~12271182). E-mail: {\tt wei229090@163.com (P.Wei), liumuhuo@163.com (M.Liu, corresponding author)}}

\vspace{3mm}

\noindent
{\large \bf Peichao Wei, \,\, Muhuo Liu$^\ast$} \\[10mm]
{\it Department of Mathematics, and Research Center for Green Development of Agriculture, South China  Agricultural  University,
Guangzhou, 510642, China} \\[3mm]

\vspace{3mm}

\noindent

\vspace{10mm}

\noindent
{\bf A B S T R A C T} \\
For a simple connected graph $G=(V,E)$, let $d(u)$ be the degree of the vertex $u$ of $G$.
The general Sombor index of $G$ is defined as
$$
SO_{\alpha}(G)=\sum_{uv\in E} \left[d(u)^2+d(v)^2\right]^\alpha
$$
where $SO(G)=SO_{0.5}(G)$ is the recently invented Sombor index.
In this paper,  we show that in the class of connected graphs with a fixed degree sequence (for which the minimum degree being equal to one),
there exists a special extremal $BFS$-graph with minimum general Sombor index
for  $0<\alpha<1$  (resp. maximum general Sombor index
for either $\alpha>1$ or $\alpha<0$).   Moreover, for any given tree, unicyclic, and bicyclic degree sequences with minimum degree 1,
there exists a unique extremal $BFS$-graph with minimum general Sombor index for $0<\alpha<1$  and  maximum general Sombor index for either $\alpha>1$ or $\alpha<0$.

\vspace{5mm}

\baselineskip=0.30in

\section{Introduction}

Throughout this paper we consider undirected simple connected graphs. Let $G$ be such a graph with vertex set
$\mathbf V(G)$ and edge set $\mathbf E(G)$.
Let $d(u) = d_G(u)$ and $N(u) = N_G(u)$ denote, respectively, the degree and neighbor set of the vertex $u \in V(G)$.
If $\mathbf V(G)=\{v_1,v_2,\ldots,v_n\}$ and $d_i=d(v_i)$, $1 \leq i \leq n$, then $\pi=(d_1,d_2,\ldots,d_n)$ is
said to be the {\bf degree sequence\/} of $G$. In what follows, we always suppose that $d_1\geq d_2\geq \cdots\geq d_n$
and denote by $\Gamma(\pi)$ the class of connected graphs with degree sequence $\pi$. A connected graph with $n$ vertices
and $n+c-1$ edges will be referred as a {\bf $c$-cyclic graph}. In particular, when $c=0$, $1$, and $2$, a $c$-cyclic
graph is also called a tree, unicyclic graph, and bicyclic graph, respectively.

Recently, Gutman proposed a geometric approach for interpreting degree-based graph invariants \cite{Gutmant1},
and according to this approach, he introduced the so-called {\bf Sombor index}, defined as,
\begin{equation}          \label{so}
SO = SO(G) = \sum_{uv\in \mathbf E(G)} \sqrt{(d(u)^2+d(v)^2)}\,.
\end{equation}
Eventually, this graph invariant attracted much attention, and in a series of researches its main mathematical
properties have been determined; see, for instance \cite{g2,g3,g4,g5,g6,g7,g9,g10,g11} and the review \cite{g1}.

One of the several modifications of the original Sombor index, Eq. (\ref{so}), is the {\bf general Sombor index\/} \cite{g12,g13},
defined as
$$
SO_\alpha = SO_\alpha(G) =\sum_{uv\in \mathbf E(G)} \big[d(u)^2+d(v)^2\big]^\alpha
$$
where $\alpha\ne 0$ is a real number. Evidently, $SO_{0.5}(G)=SO(G)$.

It is of evident interest to determine the elements of $\Gamma(\pi)$, extremal w.r.t. a certain graph invariant.
Several such researches have been published \cite{MH1,MH2,Lin1,Wang2,g14}. Among these results, in many cases the
extremal graphs are $BFS$-type ($BFS$ = breath first search).

\begin{definition}\label{21d}
Let $G$ be a connected graph. We say that $G$ is a  {\bf $BFS$-graph\/}  if there exists a vertex ordering
$v_{1}\prec v_{2}\prec \cdots\prec v_{n}$ of $\mathbf V(G)$ satisfying:
\par\noindent $(i)$  $d(v_{1})\geq d(v_{2})\geq \cdots \geq d(v_{n})$  and
$h(v_{1})\leq h(v_{2})\leq\cdots\leq h(v_{n})$, where $h(v_{i})$ is the distance between $v_i$ and $v_1$;
\par\noindent $(ii)$  let $v\in N(u)\backslash N(w)$ and $z\in N(w)\backslash N(u)$
such that $h(v)=h(u)+1$ and $h(z)=h(w)+1$. If $u\prec w$, then $v\prec z$.
\end{definition}

\begin{definition}\label{22d}{\em \cite{MH1}}
For a    $c$-cyclic   degree sequence $\pi= (d_1,d_2,. . .,d_n)$ with  $d_n = 1$ and $n\geq 3$, if $G$ is a
$BFS$-graph such that  $\{v_1, v_2, v_3\}$ forms a triangle of $G$ when $c \geq 1$, then   $G$ is called   a {\bf special extremal $BFS$-graph}.
\end{definition}

Recently, Gutman posed the problem of determining extremal graphs with minimum or maximum  Sombor index among the class of connected graphs with given degree sequence (via private communication). In this paper, we settle the minimum case. Actually, we can go further by showing the following
\begin{thm}\label{04t}
For any given degree sequence $\pi=(d_1,d_2,\ldots,d_n)$ with $d_n=1$, there exists a special extremal $BFS$-graph
with minimum $SO_{\alpha}(G)$ in the class of  $\Gamma(\pi)$ for $0<\alpha<1$, and there also  exists a special  extremal  $BFS$-graph
with maximum   $SO_{\alpha}(G)$ in the class of  $\Gamma(\pi)$ for either $\alpha>1$ or $\alpha<0$.
\end{thm}

\begin{figure}[h!]
\vspace*{0.5cm}\begin{center} \includegraphics[scale=0.50]{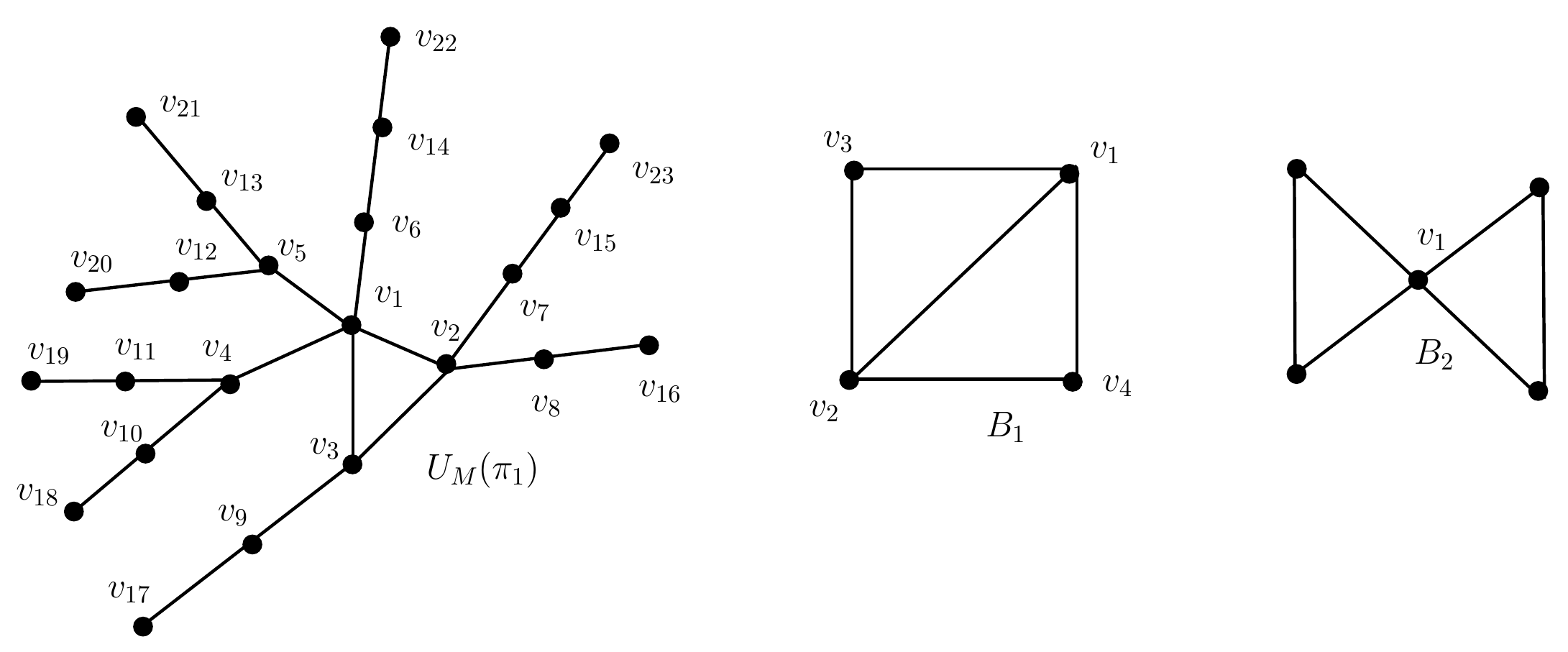} \par
\vspace*{0.0cm}
\caption{The     graphs    $U_{M}(\pi_{1})$, $B_{1}$ and $B_2$.} \label{11fig}\end{center}
     \end{figure}

Hereafter, we use the symbol $p^{(q)}$ to define $q$ copies of the real number $p$.
In \cite{Zh2008} the authors show that for any tree degree sequence $\pi$ there exists
a unique $BFS$-tree, here denoted by $T_{M}(\pi)$. The $BFS$-trees are also called {\bf greedy trees} in
the literature, e.g. \cite{Lin1,Wang2}. In fact, we can also construct a unique unicyclic $BFS$-graph $U_{M}(\pi)$
by the following breadth-first-search method for any unicyclic degree sequence $\pi=(d_{1},d_{2},\ldots,d_{n})$,
where $d_{n}=1$ \cite{MH1}: The unique cycle of $U_{M}(\pi)$ is a triangle with $V(C_{3})=\{v_{1},v_{2},v_{3}\}$.
Select the  vertex $v_{1}$  as the  root vertex and begin with $v_{1}$ of the zeroth layer. Select the vertices
$v_{2}$, $v_{3}$, $v_{4}$, $v_{5}$,\ldots, $v_{d_{1}+1}$ as the first layer such that
$N(v_{1})=\{v_{2},v_{3},v_{4},v_{5},\ldots,v_{d_{1}+1}\}$. Let   $N(v_{2})=\{v_{1},v_{3},v_{d_{1}+2},v_{d_{1}+3},\ldots,v_{d_{1}+d_{2}-1}\}$ and $N(v_{3})=\{v_{1},v_{2},v_{d_{1}+d_{2}},\ldots,v_{d_{1}+d_{2}+d_{3}-3}\}$. Then, append $d_{4}-1$ vertices to $v_{4}$ such that $N(v_{4})=\{v_{1},$$v_{d_{1}+d_{2}+d_{3}-2},\ldots,v_{d_{1}+d_{2}+d_{3}+d_{4}-4}\}$, and so on.
Informally, the $BFS$-unicyclic graph is constructed from a $BFS$-tree by adding an edge between
$v_2$ and $v_3$. As an example, considering  the unicyclic degree sequence $\pi_{1}=(5,4,3^{(3)},2^{(10)},1^{(8)})$, $U_{M}(\pi_{1})$ is depicted in Fig. \ref{11fig}.
  \par \vspace{-0.4cm}

Denote by $\mathcal {R}(G)$ the {\bf reduced graph\/} obtained  from $G$ by recursively deleting pendent vertices to the resultant
graph until no pendent vertices remain. If $c\geq 1$ and $G$ is a   $c$-cyclic graph, then $\mathcal {R}(G)$ is unique and
$\mathcal {R}(G)$ is also a   $c$-cyclic graph.
\medskip

Some paths $P_{l_{1}}$, $P_{l_{2}}$, $...$, $P_{l_{k}}$ are said to have almost equal lengths if their lengths pairwise differ at most by $1$, that is, $|l_{i}-l_{j}|\leq 1$ for $1\leq i< j\leq k$.  Let $B_1$ and $B_2$  be the two bicyclic graphs depicted in Fig. \ref{11fig}. If $\pi=(d_{1},d_{2},...,d_{n})$ is a bicyclic degree sequence with $d_n=1$, then $\sum_{i=1}^{n}d_{i}=2n+2$, which implies that $\pi$ should be one of the following two cases. When $d_n=1,$ we construct a unique  bicyclic graph $B_{M}(\pi)$ of $\Gamma(\pi)$ as follows:
\par\noindent $(i)$  If $d_{n} = 1$ and $d_{1}\geq  d_{2}\geq 3$, then let $B_{M}(\pi)$ be a $BFS$-graph such that $\mathcal {R}(B_{M}(\pi))=B_1$ and the
remaining vertices appear in a $BFS$-ordering.
\par\noindent $(ii)$  If $d_{1}\geq 5>d_{2}= 2$ and $d_{n} = 1$, then let $B_{M}(\pi)$ be the bicyclic graph with $n$ vertices obtained from $B_2$ by attaching $d_{1}-4$ paths of almost equal lengths to the $v_1$ of $B_2$ (see Fig. \ref{11fig}).

\begin{definition}\label{12d}{\em\cite{MH3}}
For a  given degree sequence $\pi=(d_1,d_2,\ldots,d_n)$ with  $d_n=1$, we say that $G$ is a
{\bf precisely  extremal graph\/} of $\Gamma(\pi)$, if $G$   has minimum  $SO_{\alpha}(G)$ among all graphs of  $\Gamma(\pi)$ for $0<\alpha<1$
and  $G$   has maximum   $SO_{\alpha}(G)$ among all graphs of  $\Gamma(\pi)$ for either $\alpha>1$ or $\alpha<0$.
\end{definition}

\begin{thm}\label{06t}
For any $c$-cyclic   degree sequence $\pi=(d_1,d_2,\ldots,d_n)$ with $d_n=1$, then  \par\noindent
$(i)$ $T_{M}(\pi)$ is a precisely extremal graph  for $c=0$;\par\noindent
$(ii)$ $U_{M}(\pi)$ is a precisely extremal   graph for $c=1$;
\par\noindent
$(iii)$ $B_{M}(\pi)$ is a precisely extremal    graph for $c=2$.
	\end{thm}

In Theorems \ref{04t} and \ref{06t}, we can only confirm that there exists a (precisely) extremal $BFS$-graph, as the   extremal  graphs of $\Gamma(\pi)$ are always not uniquely. For instance, let $\pi_2=(4,2^{(8)},1^{(4)})$ and  let $H_1$ and $H_2$ be the two trees as shown in Fig. 2. It is easily to see that $H_1$ is the unique $BFS$-tree of $\Gamma(\pi_2)$ and $SO_{\alpha}(H_1)=SO_{\alpha}(H_2)$ for any $\alpha\ne 0$.
\par\vspace*{0.8cm}
\setlength{\unitlength}{0.50mm}
\begin{picture}(40,50)
\put(105,40){\line(-1,-1){12}}
\put(93,28){\line(-1,0){30}}\put(78,28){\circle*{1.8}}\put(63,28){\circle*{1.8}}
\put(105,40){\line(1,1){12}}\put(117,52){\circle*{1.8}}
 \put(105,40){\line(-1,1){12}} \put(117,28){\circle*{1.8}}
\put(105,40){\line(1,-1){12}}
\put(93,52){\circle*{1.8}}\put(93,28){\circle*{1.8}}
\put(93,52){\line(-1,0){30}}\put(78,52){\circle*{1.8}}\put(63,52){\circle*{1.8}}
 \put(105,40){\circle*{1.8}}
\put(117,28){\line(1,0){30}}\put(132,28){\circle*{1.8}}\put(147,28){\circle*{1.8}}
\put(117,52){\line(1,0){30}}\put(132,52){\circle*{1.8}}\put(147,52){\circle*{1.8}}
 \put(100,25){$H_{1}$}

\put(205,40){\line(-1,-1){12}}
\put(193,28){\line(-1,0){15}}\put(178,28){\circle*{1.8}}
\put(205,40){\line(1,1){12}}\put(217,52){\circle*{1.8}}
 \put(205,40){\line(-1,1){12}} \put(217,28){\circle*{1.8}}
\put(205,40){\line(1,-1){12}}
\put(193,52){\circle*{1.8}}\put(193,28){\circle*{1.8}}\put(193,52){\line(-1,0){15}}\put(178,52){\circle*{1.8}}
 \put(205,40){\circle*{1.8}}
\put(217,28){\line(1,0){45}}\put(232,28){\circle*{1.8}}\put(247,28){\circle*{1.8}}\put(262,28){\circle*{1.8}}
\put(217,52){\line(1,0){45}}\put(232,52){\circle*{1.8}}\put(247,52){\circle*{1.8}}\put(262,52){\circle*{1.8}}
 \put(200,25){$H_{2}$}
 \put(80,5) { Figure 2.  The trees $H_{1}$ and $H_{2}$. }
 \end{picture}\par

The research of extremal graph in the class of connected graphs with given degree sequence
has close  relation with the Majorization theorem. Now, we introduce the notation of Majorization.

\begin{definition}{\em\cite{MA1976}}
Let $(x)=(x_{1},x_{2},...,x_{n})$ and $(y)=(y_{1},y_{2},...,y_{n})$ be two different  non-increasing
sequences of real numbers. We write $(x)\lhd (y)$ if and only if $\sum_{i=1}^{n}x_{i}=\sum_{i=1}^{n}y_{i}$,
and $\sum_{i=1}^{j}x_{i}\leq\sum_{i=1}^{j}y_{i}$ for all $j=1,2,...,n$.   The ordering $\pi\lhd \pi'$ is
said to be a {\bf majorization}.
\end{definition}

\begin{thm}\label{07t}
Let $\pi$ and $\pi'$ be two  $c$-cyclic   degree sequences and let $G$ and $G'$ be a maximum
extremal  graph  of $\Gamma(\pi)$ and $\Gamma(\pi')$, respectively. If $c\in \{0,1,2\}$ and $\pi\lhd \pi'$,
then $SO_{\alpha}(G)< SO_{\alpha}(G') $ for $\alpha>1$.
\end{thm}
\vspace{4mm}

\section{Proof of Theorems \ref{04t} and \ref{06t}}

This section is dedicated to the proofs of Theorems \ref{04t} and \ref{06t}. We need to introduce more notations.
A symmetric bivariate function $ f(x,y) $ defined on positive real numbers is called {\bf escalating\/} (resp. {\bf de-escalating\/}) if
	\begin{align}\label{21e}
		f(x_1,x_2)+f(y_1,y_2)\geq \,\,\text{(resp., $\le$ )}\,\,f(y_1,x_2)+f(x_1,y_2)
	\end{align}
holds for any $ x_1\geq y_1>0 $ and $ x_2\geq y_2>0 $, and the inequality in (\ref{21e}) is strict if $ x_1>y_1 $ and $ x_2>y_2 $.

Further, Wang \cite{Wang1}  defined the {\bf connectivity function\/} of a connected graph $G$
associated with a symmetric bivariate function $f(x,y)$ to be
$$
M_f(G)=\sum_{uv\in \mathbf E(G)}f(d(u),d(v))\,.
$$
It is easily to see that  $SO_{\alpha}(G)$  is just a  special case  of $M_f(G)$.

	The proofs of Theorems \ref{04t} and \ref{06t} rely on the following two lemmas:
\begin{lemma}\label{21l}{\em\cite{MH1}} For any  given degree sequence $\pi=(d_1,d_2,\ldots,d_n)$ with  $d_n=1$,
there exists a special extremal  $BFS$-graph $G$  such that   $M_f (G)$ is maximized   in $\Gamma(\pi)$ when $f (x, y)$ is escalating
and  $M_f(G)$ is minimized in $\Gamma(\pi)$ when $f (x,y)$ is de-escalating.
\end{lemma}
\begin{lemma}\label{22l}{\em\cite{MH1}} Let $\pi=(d_1,d_2,\ldots,d_n)$ be a   given  $c$-cyclic degree sequence
 with  $d_n=1$.   In the class of $\Gamma(\pi)$,    \par\noindent
 $(i)$  if $c=0$, then  $T_M(\pi)$ has   maximum    $M_f (G)$   when $f (x, y)$ is escalating, and
$T_M(\pi)$ has minimum  $M_f(G)$   when $f (x,y)$ is de-escalating; \par\noindent
$(ii)$  if $c=1$, then  $U_M(\pi)$ has   maximum    $M_f (G)$   when $f (x, y)$ is escalating, and
$U_M(\pi)$ has minimum  $M_f(G)$   when $f (x,y)$ is de-escalating; \par\noindent
$(iii)$  if $c=2$, then  $B_M(\pi)$ has   maximum    $M_f (G)$   when $f (x, y)$ is escalating, and
$B_M(\pi)$ has minimum  $M_f(G)$   when $f (x,y)$ is de-escalating;
\end{lemma}
	
Throughout this  paper, denote  by $h(x,y)=\big(x^2+y^2\big)^{\alpha}$, where $\min\{x,y\}>0$.
From Lemmas \ref{21l}--\ref{22l}, to show Theorems \ref{04t} and \ref{06t},
it suffices to show that the following  proposition holds.

\begin{proposition}\label{21p}    $SO_{\alpha}(G)$ is escalating for  $\alpha>1$ or $\alpha<0$, and  $SO_{\alpha}(G)$ is de-escalating for  $0<\alpha<1$.
\end{proposition}
\noindent{\bf Proof.} In what follows, we suppose that $x_1\ge y_1\ge 1$ and $x_2\ge y_2\ge 0$.  It suffices to show that $h(x,y)$ satisfies  \eqref{21e}. Since the equality holds in  \eqref{21e} for  either $x_1=y_1$ or  $x_2=y_2$, we may suppose that $x_1>y_1$ and $x_2>y_2$.

 One can easily check that \begin{align}\nonumber  &h(x_1,x_2)+h(y_1,y_2)-h(y_1,x_2)-h(x_1,y_2)\\[2mm] \label{31e}=&\int_{y_1}^{x_1} 2t\alpha(t^2+x_2^2)^{\alpha-1} dt-\int_{y_1}^{x_1} 2t\alpha(t^2+y_2^2)^{\alpha-1} dt.\end{align}
 \par\noindent {\bf Case 1. $0<\alpha <1$.} Since $2t\alpha(t^2+x_2^2)^{\alpha-1}>0$ and $2t\alpha(t^2+y_2^2)^{\alpha-1}>0$ for $t\ge y_1>0$, we have  \begin{align}\label{32e}2t\alpha(t^2+x_2^2)^{\alpha-1}<2t\alpha(t^2+y_2^2)^{\alpha-1},\end{align} as $x_2>y_2>0$ and $0<\alpha<1$. Combining \eqref{32e} with $x_1>y_1>0$, we can conclude that $h(x_1,x_2)+h(y_1,y_2)<h(y_1,x_2)+h(x_1,y_2)$ by \eqref{31e}. Thus, $SO_{\alpha}(G)$ is de-escalating for  $0<\alpha<1$.   \par\medskip\noindent
 {\bf Case 2. $\alpha <0$ or $\alpha >1$.}  If $\alpha>1$, then $2t\alpha(t^2+x_2^2)^{\alpha-1}>2t\alpha(t^2+y_2^2)^{\alpha-1}$ for
 $t\ge y_1>0$ and $x_2>y_2$. Combining this with   $x_1>y_1>0$, we can conclude that $h(x_1,x_2)+h(y_1,y_2)>h(y_1,x_2)+h(x_1,y_2)$ by \eqref{31e}, which implies that $SO_{\alpha}(G)$ is escalating for  $\alpha>1$.

 Otherwise,   $\alpha<0$. Since  $2t\alpha(t^2+x_2^2)^{\alpha-1}<0$ and $2t\alpha(t^2+y_2^2)^{\alpha-1}<0$, we have $2t\alpha(t^2+x_2^2)^{\alpha-1}>2t\alpha(t^2+y_2^2)^{\alpha-1}$ for $t\ge y_1>0$ and $x_2>y_2>0$. Taking this with $x_1>y_1>0$ into consideration, we also deduce that  $SO_{\alpha}(G)$ is escalating for  $\alpha<0.$
	\qed

\section{Proof of Theorem \ref{07t}}
 The following definition  will play a crucial role in the proof of Theorem \ref{07t}.

 \begin{definition}{\em\cite{MH1}}
 A non-negative escalating function $ f(x,y) $ is called a
 {\bf good escalating} function, if $ f(x,y) $ satisfies  $\frac{\partial f(x,y)}{\partial x}>0$,
 $\frac{\partial^2 f(x,y)}{\partial x^2}\ge 0$, and  $$f(x_1+1,x_2)+f(x_1+1,y_2)+f(x_1+1,y_1-1)>f(x_1,x_2)+f(y_1,y_2)+f(x_1,y_1)$$
 holds for any $x_1\ge y_1\ge 2$ and $x_2\ge y_2\ge 1$.
 \end{definition}

\begin{lemma}\label{23l}{\em\cite{MH2}}
Let $\pi$ and $\pi'$ be two  $c$-cyclic   degree sequences with $\pi\lhd \pi'$ and $c\in \{0,1,2\}$. Let
 $G$ and $G'$ have maximum $M_f(G)$ and $M_f(G')$ in the class of $\Gamma(\pi)$ and $\Gamma(\pi')$, respectively.
 If $f (x, y)$ is good  escalating, then  $M_f(G)<M_f(G')$.
\end{lemma}	

\noindent{\bf Proof of Theorem \ref{07t}:} By Lemma \ref{23l}, to complete the proof of Theorem \ref{07t}, it suffices to show that
	 $h(x,y)$ is a good escalating function for  $\alpha>1$.  We already know that $ h(x,y) $ is a non-negative escalating function for $ \alpha>1 $ by Proposition \ref{21p}.    \par

 Since  $\alpha>1$, we have  $$\text{ $\frac{\partial h(x,y)}{\partial x}=2x\alpha (x^2+y^2)^{\alpha -1}>0$ and  $\frac{\partial^2 h(x,y)}{\partial x^2}=2\alpha (x^2+y^2)^{\alpha -2}\big[(x^2+y^2)+2x^2(\alpha -1)\big] {>} 0$.}$$  Since  $h(x,y)$ is a strictly  increasing function  on   $x$, we have  $h(x_1+1,x_2)>h(x_1,x_2)$ and  $h(x_1+1,y_2)>h(y_1,y_2)$, as $x_1\ge y_1$.

 Since $\big[(x_1+1)^2+(y_1-1)^2\big]-\big[x_1^2+y_1^2\big]=2(x_1-y_1+1)>0$  for $x_1\ge y_1$, we have    $$h(x_1+1,y_1-1)-h(x_1,y_1)=\big[(x_1+1)^2+(y_1-1)^2\big]^\alpha -\big(x_1^2+y_1^2\big)^\alpha>0$$  for $\alpha>1$.  Now, we can see that $h(x,y)$ is a good escalating function for  $\alpha>1$.  \qed
\par\medskip\noindent
\noindent
 {\it Acknowledgement\/}.  The authors would like  to thank  Professor Ivan Gutman for guiding us to the research of the Sombor index.

\end{document}